\documentclass[a4paper,12pt]{article}
\usepackage[latin1]{inputenc}
\usepackage{graphicx}
\usepackage{amsmath,amssymb}
\usepackage{amsfonts}

\newtheorem{thm}{Theorem}

\newtheorem{la}[thm]{Lemma}

\newtheorem{conj}{Conjecture}

\DeclareMathOperator{\uv}{\lower0.8ex\hbox{\includegraphics[height=4mm]{arc_uv.eps}}}

\newcommand{\qedgo}{\vspace{0.35cm}}

\newcommand{\tume}{\hspace*{0.2em}}

\newenvironment{u-enumerate}{%
  \begin{enumerate}%
 }{\end{enumerate}%
}

\newenvironment{Proofof}[1]{\expandafter
\paragraph{\rm{\it Proof of\/} #1.}}{\hfill $\blacksquare$\qedgo}

\begin{document}
\begin{center}
{\Large On the balanced decomposition number}
\end{center}
\begin{center}
{\large Tadashi Sakuma\footnote{Email: {\texttt{sakuma@e.yamagata-u.ac.jp}}}
}
\end{center}
\begin{center}
{
Systems Science and Information Studies \\%
Faculty of Education, Art and Science \\%
Yamagata University \\%
1-4-12 Kojirakawa, Yamagata 990-8560, Japan
}
\end{center}
\begin{abstract}
A {\em balanced coloring} of a graph $G$ means a triple 
$\{P_1,P_2,X\}$ of mutually disjoint subsets of the vertex-set $V(G)$ 
such that $V(G)=P_1 \uplus P_2 \uplus X$ and $|P_1|=|P_2|$. 
A {\em balanced decomposition} associated with 
the balanced coloring $V(G)=P_1 \uplus P_2 \uplus X$ of $G$ is 
defined as a partition of $V(G)=V_1 \uplus \cdots \uplus V_r$  
(for some $r$) such that, for every $i \in \{1,\cdots,r\}$, the subgraph 
$G[V_i]$ of $G$ is connected and $|V_i \cap P_1| = |V_i \cap P_2|$.  
Then the {\em balanced decomposition number} of a graph 
$G$ is defined as the minimum integer $s$ such that, for every balanced 
coloring $V(G)=P_1 \uplus P_2 \uplus X$ of $G$, there exists a balanced 
decomposition $V(G)=V_1 \uplus \cdots \uplus V_r$ whose every 
element $V_i (i=1, \cdots, r)$ has at most $s$ vertices. 
S. Fujita and H. Liu [\/SIAM J. Discrete Math. 24, (2010), pp. 
1597--1616\/] proved a nice theorem which states that the 
balanced decomposition number of a graph $G$ is at most $3$ 
if and only if $G$ is $\lfloor\frac{|V(G)|}{2}\rfloor$-connected. 
Unfortunately, their proof is lengthy (about 10 pages) and complicated. 
Here we give an immediate proof of the theorem. This proof 
makes clear a relationship between balanced decomposition 
number and graph matching. 
\end{abstract}
keywords: graph decomposition, coloring, connectivity, bipartite matching

\section{Introduction}

Throughout this paper, we only consider finite undirected graphs with no 
multiple edges or loops. For a graph $G$, let $V(G)$ and $E(G)$ denote
the vertex-set of $G$ and the edge-set of $G$, respectively. For a subset
$X \subseteq V(G)$,  $G[X]$ denotes the subgraph of $G$ induced by $X$, 
and $N_G(X)$ denotes the set $\big\{y \in V(G) \setminus X | \exists x \in 
X, \{x,y\} \in E(G)\big\}$. This set $N_G(X)$ is called the 
{\em open neighborhood} of $X$ in $G$. 
A subset $Y \subseteq V(G)$ is called a {\em vertex-cut} of $G$ if 
there is a partition $V(G)\setminus Y = X_1 \uplus X_2$ such that 
$|X_i| \geqq 1$ and $N_{G[V(G) \setminus Y]}(X_i) =\emptyset ~ (i=1,2)$. 
For other basic definitions in graph theory, please consult \cite{D}. 

In 2008, S. Fujita and T. Nakamigawa\tume\cite{FN} introduced a new 
graph invariant, namely the {\em balanced decomposition number} 
of a graph, which was motivated by the estimation of the number of 
steps for pebble motion on graphs. 
A {\em balanced coloring} of a graph $G$ means a triple 
$\{P_1,P_2,X\}$ of mutually disjoint subsets of $V(G)$ such that 
$V(G)=P_1 \uplus P_2 \uplus X$ and $|P_1|=|P_2|$. 
Then a {\em balanced decomposition} of $G$ 
associated with its balanced coloring $V(G)=P_1 \uplus P_2 \uplus X$ 
is defined as a partition of $V(G)=V_1 \uplus \cdots \uplus V_r$  (for 
some $r$) such that, for every $i \in \{1,\cdots,r\}$, $G[V_i]$ is
connected and $|V_i \cap P_1| = |V_i \cap P_2|$.  Note that every 
disconnected graph has a balanced coloring which admits no 
balanced decompositions. Now the {\em balanced 
decomposition number} of a connected graph $G$ is defined as 
the minimum integer $s$ such that, for every  balanced coloring 
$V(G)=P_1 \uplus P_2 \uplus X$ of $G$, there exists a balanced 
decomposition $V(G)=V_1 \uplus \cdots \uplus V_r$ whose every 
element $V_i (i=1, \cdots, r)$ has at most $s$ vertices.

The set of the starting and the target arrangements of mutually 
indistinguishable pebbles on a graph $G$ can be modeled as a 
balanced coloring $V(G)=P_1 \uplus P_2 \uplus X$ of $G$. Then, 
as is pointed out in \cite{FN}, the balanced decomposition number 
of $G$ gives us an upper-bound for the minimum number of 
necessary steps to the pebble motion problem, and, for several 
graph-classes, this upper bound is sharp. 

In addition to the initial motivations and their applications in 
\cite{FN}, this newcomer graph invariant turns out to have 
deep connections to some essential graph theoretical concepts. 
For example, the following conjecture in \cite{FN} indicates a 
relationship between this invariant and the vertex-connectivity 
of graphs: 

\begin{conj}[S. Fujita and T. Nakamigawa (2008)]
The balanced decomposition number of $G$ is at most
$\lfloor\frac{|V(G)|}{2}\rfloor + 1$ if $G$ is 
$2$-connected. 
\end{conj}
Recently, G. J. Chang and N. Narayanan\tume\cite{CN} 
announced a solution to this conjecture. 

Then especially, S. Fujita and H. Liu\tume\cite{FL} proved 
the affirmation of the ``high''-connectivity counterpart of 
the above conjecture, as follows: 

\begin{thm}[S. Fujita and H. Liu (2010)]
Let $G$ be a connected graph with at least $3$ vertices. Then 
the balanced decomposition number of $G$ is at most $3$ 
if and only if $G$ is $\lfloor\frac{|V(G)|}{2}\rfloor$-connected. 
\end{thm}\label{thmFL}

Thus, there may be a trade-off between the vertex-connectivity and 
the balanced decomposition number. This interesting relationship 
should be investigated for its own sake. 

Unfortunately, the proof of Theorem\tume\ref{thmFL} 
in \cite{FL} is lengthy (about $10$ pages) and complicated. 

In this note, we give a new proof of the theorem\tume\ref{thmFL}. 
The advantages of our proof is that it is immediate and makes clear 
a relationship between balanced decomposition number and graph 
matching.

\section{A quick proof of Theorem\tume\ref{thmFL}}

We show our proof of the theorem\tume\ref{thmFL} here. 

\begin{Proofof}{Theorem\tume\ref{thmFL}}
In order to prove the {\bf if part}, let us define the following
new bipartite graph $H$ from a given balanced coloring $V(G)=P_1 
\uplus P_2 \uplus X$ of a graph $G$: 
\begin{enumerate}
      \item The partite sets of $H$ are $V_1(H):=P_1 \uplus X_1$ 
                  and $V_2(H):=P_2 \uplus X_2$, where each $X_i:=\{(x,i) ~ |
                  ~ x \in X\} ~ (i=1,2)$ is a copy of the set $X(\subseteq V(G))$. 
      \item The edge set $E(H)$ of $H$ is defined as follows: \vspace*{-3pt}
                 \begin{equation*}
                  \begin{split}
                   &E(H):=\big\{\{p_1,p_2\} ~|~ p_1 \in P_1, p_2 \in P_2,
                   \{p_1,p_2\}\in E(G)\big\} \\
                   &\qquad \qquad \cup\big\{\{p_1,(x,2)\}\ |~ p_1 \in P_1, x \in X,
                   \{p_1, x\} \in E(G) \big\} \\
                   &\qquad \qquad \cup\big\{\{(x,1),p_2\}\ |~ x \in X, p_2 \in P_2,
                   \{x, p_2\} \in E(G) \big\} \\
                   &\qquad \qquad \cup\big\{\{(x,1),(x,2)\}\ |~ x \in X \big\}.
                  \end{split}
                 \end{equation*}
\end{enumerate}
Then clearly, the balanced coloring $V(G)=P_1 \uplus P_2 \uplus X$ of $G$ has 
a balanced decomposition $V(G)=V_1 \uplus \cdots \uplus V_r$  whose 
every element $V_i (i=1, \ldots, r)$ consists of at most $3$ vertices, 
if and only if the graph $H$ has a perfect matching. 
Then we use here the famous ``Hall's Marriage Theorem''\cite{H}, as follows. 
\begin{la}[P. Hall(1935)]
Let $G$ be a bipartite graph whose partite sets are $V_1(G)$ and
$V_2(G)$. Suppose that $|V_1(G)|=|V_2(G)|$. Then $G$ has a perfect
matching if and only if every subset $U$ of $V_1(G)$ satisfies 
$|U| \leqq |N_G(U)|$. 
\end{la}\label{Hall}
Now, suppose that $H$ does not have any perfect matching. Then, from 
lemma\tume\ref{Hall}, $\exists A \subseteq P_1,
\exists B \subseteq X_1, |N_H(A \cup B)| \leqq |A| + |B| - 1$. Let $C:=P_2
\setminus N_H(A \cup B)$ and $D:=X_2 \setminus N_H(A \cup B)$. Then, 
by symmetry, $|N_H(C \cup D)| \leqq |C| + |D| - 1$ also holds. Furthermore, 
by the definition of $H$, $|B| \leqq |X_2 \setminus D|$ and 
$|D| \leqq |X_1 \setminus B|$ hold, and hence $0 \leqq |X| - |B| -|D| \leqq  
|A| + |C| - |P_1| -1 = |A| + |C| - |P_2| -1$ satisfies. Please see 
Figure\tume\ref{fig:Hall} which shows this situation. 
\begin{figure}[htbp]
   \begin{center}
      \rotatebox{-90}{\includegraphics[width=9cm]{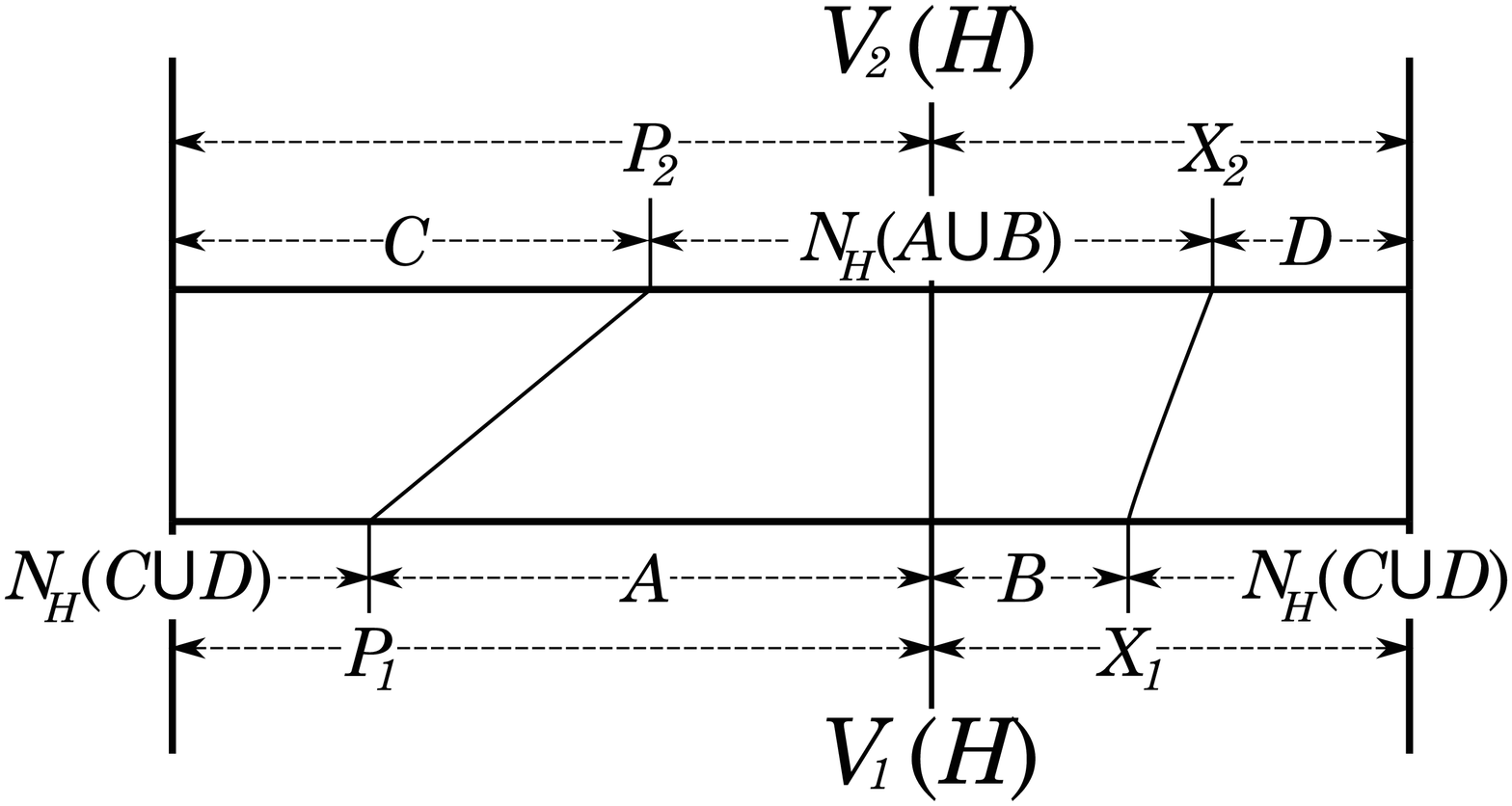}}
      \caption{The bipartite graph $H$ which has no perfect matching.}
     \label{fig:Hall}
   \end{center}
\end{figure}
The vertex-cut of $V(G)$
corresponding to the set $(P_1 \setminus A) \cup (P_2 \setminus C) \cup
(X_1 \setminus B)$ separates $G[C]$ from its remainder. 
By symmetry, the vertex-cut of $V(G)$ corresponding to the set 
$(P_1 \setminus A) \cup (P_2 \setminus C) \cup
(X_2 \setminus D)$ separates $G[A]$ from its remainder. 
Hence if $G$ is  $\lfloor\frac{|V(G)|}{2}\rfloor$-connected, 
$|V(G)| -1 \leqq 2(|P_1| - |A| + |P_2| - |C|) + (|X| - |B|) + (|X| - |D|) = 
(|P_1| + |P_2| + |X|)  - 2\big((|A| + |C| - |P_1|) - (|X| - |B| - |D|)\big) -
(|X| - |B| - |D|) \leqq |V(G)| - 2$, a contradiction. 

The proof of the {\bf only if part} is given by a construction of 
special balanced colorings, which is the same as the original 
one in \cite{FL}. We will transcribe the construction 
only for the convenience of readers. 

Suppose that $G$ is not $\lfloor\frac{|V(G)|}{2}\rfloor$-connected. 
And let $Y$ denote a minimum vertex-cut of $G$. Note that 
$2|Y| \leqq |V(G)| -2$. Then $G[V(G) \setminus Y]$ is divided into two 
graphs $G_1$ and $G_2$ such that $|V(G_i)| \geqq 1$ and 
$N_{G[V(G) \setminus Y]}(V(G_i))=\emptyset ~(i=1,2)$. 
Without loss of generality, we assume that 
$|V(G_1)| \leqq |V(G_2)|$. Let $l$ denote the number 
$\min\{|Y|, |V(G_1)|-1\}$. Suppose an arbitrary balanced coloring 
$V(G)=P_1 \uplus P_2 \uplus X$ of $G$ such that $|Y \cap P_1|=l$ and 
$|Y \cap P_2| = |Y| - l$ and $|V(G_1) \cap P_2| = l+1$ and $V(G_1) \cap
P_1 = \emptyset$. Then, it is easy to see that every balanced decomposition 
associated with such a balanced coloring has at least one component whose 
vertex-size is at least $4$, that is, the balanced decomposition 
number of $G$ is at least $4$. 
\end{Proofof}


%
\end{document}